\documentclass[+2pt]{amsart}

\usepackage{amsmath,amssymb}
\usepackage{mathpazo}
\usepackage{amsthm}
\usepackage{comment}
\usepackage{amscd}
\usepackage[all]{xy}
\usepackage{cite}

\usepackage{graphicx}

\usepackage{eq}

\theoremstyle{plain}
\newtheorem{thm}{Theorem}[section]
\newtheorem{lem}[thm]{Lemma}
\newtheorem{prop}[thm]{Proposition}
\newtheorem{cor}[thm]{Corollay}

\theoremstyle{definition}
\newtheorem{mydef}[thm]{Definition}

\newcommand{\Oct}{\mathbb{O}}

\newcommand{\Aff}{\mathrm{Aff}}

\newcommand{\J}{\mathcal{J}}
\newcommand{\GL}{\mathrm{GL}}

\newcommand{\minimat}{( \begin{smallmatrix} a&b\\ c&d
\end{smallmatrix} \bigr)}

\newcommand{\oline}[1]{\overline{#1}}

\newcommand{\ks}{k^{\mathrm{sep}}}

\newcommand{\at}{\mathfrak{t}}

\newcommand{\Mj}{\mathcal{M}} 
\newcommand{\an}{\mathfrak{n}}

\newcommand{\Trace}{\mathrm{Tr}}
\begin{document}

\address[R. Kato]
{Department of Mathematics, Graduate School of Science, 
Kyoto University, Kyoto 606-8502, Japan}
\email{rkato@math.kyoto-u.ac.jp}

\address[A. Yukie]
{Department of Mathematics, Graduate School of Science, 
Kyoto University, Kyoto 606-8502, Japan}
\email{yukie@math.kyoto-u.ac.jp}
\thanks{The first author was partially supported by 
Grant-in-Aid (B) (24340001) \\}


\title{On equivariant maps related to the space of 
pairs of exceptional Jordan algebras}
\author{Ryo Kato}
\author{Akihiko Yukie}

\maketitle


\begin{abstract}
Let $\J$ be the exceptional Jordan algebra and $V=\J\oplus \J$.
We construct an equivariant map from $V$ to $\mathrm{Hom}_k(\J\otimes
\J,\J)$ defined by homogeneous polynomials of degree $8$ such that if
$x\in V$ is a generic point, then the image of $x$
is the structure constant of the isotope of $\J$ corresponding to $x$.
We also give an alternative way to define the isotope corresponding to
a generic point of $\J$ by an equivariant map from $\J$ to the space
of trilinear forms.
\end{abstract}

\section{Introduction}
\label{sec:introduction}

Let $k$ be a field of characteristic not equal to $2, 3$, $\ks$ 
the separable closure of $k$ and $\oline{k}$ 
the algebraic closure of $k$.
Let $\widetilde\Oct$ be the split octonion over $k$.  
It is the normed algebra
over $k$ obtained by the Cayley--Dickson process
(see \cite[pp.101--110]{harvey}). If $A$ is the algebra
of $2\times 2$ matrices, $\widetilde\Oct$ is $A(+)$ with the notation 
of \cite{harvey}.  An {\it octonion} is, by definition, a normed algebra
which is a $k$-form of $\widetilde\Oct$.
Let $\Oct$ be an octonion. We use the notation $\|x\|$ 
for the norm of $x\in \Oct$.  If $a\in k$, $\|ax\|=a^2\|x\|$. 
Also if $x,y\in \Oct$, then $\|xy\|=\|x\| \|y\|$. 
For $x, y \in \Oct$, let
\begin{equation*}
Q(x,y) = \tfrac 12(\|x+y\|-\|x\|-\|y\|).  
\end{equation*}
This is a non-degenerate symmetric bilinear form such that $Q(x,x)=\|x\|$. 
Let $W\subset \Oct$ be the orthogonal complement of $k\cdot 1$
with respect to $Q$.  If $x=x_1+x_2$  where $x_1\in k\cdot 1$,
$x_2\in W$, then we define $\oline{x}=x_1-x_2$ and call it the
{\it conjugate} of $x$. Note that $\|x\|=x\oline{x}$. 
For $x \in \Oct$, we define the trace tr$(x)$ by tr$(x) = x + \oline{x}$.
It is easy to verify that 
\begin{equation*}
\mathrm{tr}(xy) =
\mathrm{tr}(yx),\;
2Q(x, y) = \mathrm{tr}(x\oline y). 
\end{equation*}

Let $\GL(n)$ be the group of $n\times n$ invertible matrices 
over $k$. If $V$ is a finite dimensional vector space, then 
we denote the group of invertibel linear maps from $V$ to $V$
by $\GL(V)$. 
Let  $\J$ be the exceptional Jordan algebra over $k$.
Elements of $\J$ are of the form:
\begin{equation*}
X = 
\begin{pmatrix}
s_1 & x_3& \oline{x_2}\\
\oline{x_3} & s_2 & x_1\\
x_2 & \oline{x_1} & s_3
\end{pmatrix},
\ \ s_i \in  k, \ x_i \in \Oct \ \ (i =1, 2, 3),
\end{equation*}
The multiplication  of $\J$ is defined as follows:
\begin{equation*}
X \circ Y = \frac{1}{2} (XY + YX), 
\end{equation*}
where the multiplication used on the right-hand side is 
the multiplication of matrices.

The algebraic groups $E_6$ and $GE_6$ are given by
\begin{align*}
E_6 & = \{ L \in  \mathrm{GL}(\J)  |  {}^{\forall}
 X \in \J, \  \det(LX) = \det(X) \}, \\
GE_6 & = \{ L \in  \mathrm{GL}(\J)  | {}^{\forall}
 X \in \J,   \det(LX) = c(L)\det(X) \
\text{for some} \  c(L) \in \GL (1)  \}
\end{align*}
respectively.
Then $c: GE_6 \to \GL (1)$ is a character and there exists an exact sequence
\begin{equation}\label{E6exact}
0 \to E_6 \hookrightarrow GE_6 \overset{c}{\to} \GL (1) \to 0. 
\end{equation}

It is known that $E_6$ is a smooth connected quasi-simple
simply-connected algebraic group of type ${\mathrm E}_6$ (see \cite[p.181,
Theorem 7.3.2]{Springer}).
The terminology ``quasi-simple'' means that
its inner automorphism group is simple (see
\cite[p.136]{Springer-LAG}). 
The smoothness of the group follows from the fact that 
the dimension of $E_6$ as a variety and the dimension of the 
Lie algebra of $E_6$ coincide (see the proof of \cite[p.181,
Theorem 7.3.2]{Springer}).

Let $H_1=E_6$, $G_1 = GE_6$,  
$H = H_1 \times \GL(2)$ and $G = G_1 \times \GL(2)$. 
Let $V=\J\otimes \Aff^2$.  
We regard elements of $V$ as the set of $x =x_1 v_1 + x_2 v_2$ where
$x_1, x_2 \in \J$ and two variables $v_1, v_2$.
The action of $g = (g_1, g_2) \in G$ where $g_2 = \minimat$ on $V$ is 
given by 
\begin{equation}
g(x_1 v_1 + x_2 v_2) = g_1(x_1)(av_1 + cv_2) +  g_1(x_2)(bv_1 + dv_2).
\end{equation}
For $x=x_1v_1+x_2v_2\in V$, we put $F_x(v)=F_x(v_1,v_2) = \det(x)$. 
Then $F_x(v)$ is a binary cubic form. Let $\Delta(x)$ be the 
discriminant of $F_x(v)$ as a polynomial of $v$. 
Let
\begin{equation}
\label{eq:w-defn}
w = w_1v_1 + w_2 v_2 =
\begin{pmatrix}
1&0&0\\
0&-1&0\\
0&0&0
\end{pmatrix}
v_1 +
\begin{pmatrix}
0&0&0\\
0&1&0\\
0&0&-1
\end{pmatrix}
v_2
\ \in \J \otimes \mathrm{Aff}^2.
\end{equation}
Then it is easy to see that $F_w(v)=v_1v_2(v_1-v_2)$ and 
$\Delta(w)=1$. 
The pair $(G,V)$ is an irreducible regular 
{\it prehomogenerous vector space}
(see \cite[Proposition 3.2]{kato-yukie-jordan} and \cite{regularity-of-PVS})
and the polynimial $\Delta(x)$ is what we call a {\it relative invariant 
polynomial}.  We define $V^{\mathrm{ss}}=\{x\in V\ | \ \Delta(x)\not=0 \}$.
Points in $V^{\mathrm{ss}}$ are called {\it semi-stable points}. 
As we pointed out as above, $w\in  V^{\mathrm{ss}}_k$. 
The polynomial $\Delta(x)$ is of degree $12$. 
If we put $\chi(g)=c(g_1)^4\det(g_2)^6$ for $g=(g_1,g_2)\in G$, 
then $\Delta(gx)=\chi(g)\Delta(x)$.

A Jordan algebra $\Mj$ is called an isotope of $\J$
if $\Mj\otimes k^{\mathrm{sep}}\cong \J\otimes k^{\mathrm{sep}}$
and the ``determinant'' of $\Mj$ is a constant multiple
of that of $\J$. For details, the reader should see 
\cite[pp.154--158]{Springer}.
If $\Mj_1,\Mj_2$ are isotopes of $\J$ and
$\an_1\subset \Mj_1,\an_2\subset \Mj_2$ are
cubic \'etale subalgebras, then the pairs
$(\Mj_1,\an_1)$, $(\Mj_2,\an_2)$ are defined to
be equivalent if there exists a $k$-isomorphism
$\Mj_1\to \Mj_2$ which induces an isomorphism
from $\an_1$ to $\an_2$. 
Let $\mathrm{JIC}(k)$ be the set of equivalence classes
of pairs $(\Mj,\an)$ as above.

In \cite{kato-yukie-jordan}, 
to each point in $V^{\mathrm{ss}}_k$,
a pair $(\Mj,\an)\in \mathrm{JIC}(k)$ 
was associated.
It is proved in \cite[Theorem 5.8]{kato-yukie-jordan} that 
there is a bijective correspondence between the set  
$G_k\backslash V^{\mathrm{ss}}_k$ of rational orbits 
and $\mathrm{JIC}(k)$. Moreover, an equivariant map
$m:V\to \J$ was defined and the Jordan algebra corresponding
to $x\in  V^{\mathrm{ss}}_k$ was explicitly described using
the point $m(x)$ (see \cite[Section 5]{kato-yukie-jordan}).

Let $\mathfrak t\subset \J$ be the
subalgebra of diagonal matrices, which is
isomorphic to $k^3$.
Let $\mathrm{Aut}(\J)$ be the
algebraic group of automorphisms of
the Jordan algebra $\J$. We define 
$\mathrm{Aut}(\J,\mathfrak t)$ 
to be the subgroup of
$\mathrm{Aut}(\J)$ consisting of automorphisms
$L$ such that $L(\mathfrak t)=\mathfrak t$. 
It is proved in (\cite[Theorem 3.1, Lemma 4.2]{kato-yukie-jordan})
that $G_w\cong \GL(1)\times \mathrm{Aut}(\J,\mathfrak t)$. 
The first Galois cohomology set 
\begin{math}
{\mathrm H}^1(k,\GL(1)\times \mathrm{Aut}(\J,\mathfrak t))
\end{math}
can be identified with 
\begin{math}
{\mathrm H}^1(k,\mathrm{Aut}(\J,\mathfrak t)) 
\end{math}
and there is a natural map
\begin{math}
{\mathrm H}^1(k,\mathrm{Aut}(\J,\mathfrak t)) 
\to {\mathrm H}^1(k,\mathrm{Aut}(\J)). 
\end{math}
One can show by standard argument that elements of
${\mathrm H}^1(k,\mathrm{Aut}(\J))$
correspond bijectively with 
$k$-forms of $\J$.

Suppose that $x = g_xw\in  V^{\mathrm{ss}}_k$
where $g_x\in G_{k^{\mathrm{sep}}}$.
Then
\begin{math}
h_x:\mathrm{Gal}(k^{\mathrm{sep}}/k)\ni\sigma
\mapsto g_x^{-1}g_x^{\sigma}\in G_{w\, k^{\mathrm{sep}}}
\end{math}
is a 1-cocycle, which defines an element, say $c_x$
of ${\mathrm H}^1(k,G_w)$.  This element $c_x$ 
does not depend on the choice of $g_x$.
Since there is a natural map
${\mathrm H}^1(k,G_w)\cong {\mathrm H}^1(k,{\mathrm{Aut}}(\J,\at))\to
 {\mathrm H}^1(k,{\mathrm{Aut}}(\J))$ ($\cong$ means bijection), 
$c_x$ determines a $k$-form of $\J$. 
If $x$ corresponds to a pair $(\Mj,\an)$,
then $\Mj$ is the $k$-form of $\J$
which is determined by $c_x$. 
The underlying vector space of $\Mj$ is
$\J$. To define a Jordan algebra sutructure
on $\J$, it is enough to define the product
strucuture, which is 
given by an element  of
$\mathrm{Hom}_k(\J\otimes \J,\J)$
(we call this element the ``structure constant'').

If we choose bases of $V,\J$ as $k$-vector spaces, 
the map $m:V\to \J$ is given by
homogeneous polynomials of degree $4$ on $V$. 
The structure constant which is associated to elements of
$\J$ is, with the denominator multiplied,  
given by homogeneous polynomials of degree $11$ on $\J$.  
So the structure constant of $\Mj$ is given by
homogeneous polynomials of degree $44$ on $V$.

The first purpose of this paper is to prove the
following theorem (see Section \ref{sec:equivariant-mapI}).

\begin{thm}
\label{thm:main-theorem-equiv}
There is an equivariant map
\begin{math}
\mathrm{S}:V\ni x \mapsto \mathrm{S}_x\in \mathrm{Hom}_k(\J\otimes \J,\J)
\end{math}
defined by homogeneous polynomials of degree $8$ on $V$
such that if $x \in  V^{\mathrm{ss}}_k$
corresponds to the pair $(\Mj,\an)\in {\mathrm{JIC}}(k)$,  then 
$\Delta(x)^{-1}\mathrm{S}_x$ is the structure constant of  
the Jordan algebra $\Mj$. 
\end{thm}

Let $a\in\J$ and $\det(a)\not=0$. 
In \cite[p. 155]{Springer}, for $X,Y\in \J$, 
the product structure $X\circ_a Y$ defining  
the isotope $\J_a$ and 
the corresponding symmeric bilinear form 
$Q_a(X,Y)$ are given (see (\ref{eq:bilinear-product})).  
Then 
$\J^3\ni (X,Y,Z)\mapsto {\mathrm{T}}_a(X,Y,Z) = Q_a(X\circ_a Y,Z)$
is a trilinear form on $\J$. 
To construct $T_a$ first equivariantly 
to provide an alternative way to define 
the product structure on $\J_a$ is another purpose
of this paper. 

We prove the following theorem
in Section \ref{sec:equivariant-mapII}.

\begin{thm}
\label{thm:main-theorem-equiv2}
There is an equivariant map
\begin{math}
\mathrm{T}:\J\ni a \mapsto \mathrm{T}_a\in \mathrm{Hom}_k(\J\otimes \J\otimes \J,k)
\end{math}
defined by homogeneous olnynomials of degree $6$ 
such that if for $X,Y\in\J$, $X\circ_a Y\in\J$ is the element such 
$Q_a(X\circ_a Y,Z)=\det(a)^{-1}{\mathrm T}_a(X,Y,Z)$ for all $Z\in\J$, 
then this product structure coincides with that of the isotope 
$\J_a$. 
\end{thm}


\newcommand{\lan}{\langle}
\newcommand{\ran}{\rangle}
\newcommand{\ccd}{,\ldots,}

\section{Equivariant map I}
\label{sec:equivariant-mapI}

We prove Theorem \ref{thm:main-theorem-equiv} in this section. 
We first define basic notions and then define 
the desired equivariant map. 

We denote the $3\times 3$ diagonal matrix 
with diagonal entries $\alpha_1,\alpha_2,\alpha_3\in k^{\times}$
by $\mathrm{diag}(\alpha_1,\alpha_2,\alpha_3)$. 
$\J^{n\otimes}$ is the tensor product of $n$ copies of $\J$. 
We define a symmetric trilinear form 
$D$ on $\J$ by
\begin{equation*}
\begin{split}
6D(X, Y , Z)  = 
&\det(X + Y  + Z) - \det(X + Y) - \det(Y  + Z) - \det(Z + X) \\
&+ \det(X) + \det(Y) +\det(Z)
\end{split}
\end{equation*}
for $X,Y,Z\in\J$. 
Let $\Trace(X)$ be the sum of diagonal entries of $X\in\J$.  
For $X, Y \in \J$, 
we define a symmetric bilinear form $\langle \ ,\ \rangle $ on $\J$ by
\begin{equation*}
\langle  X, Y\rangle  = \Trace(X\circ Y).
\end{equation*}
One can verify by direct computation that the symmetric bilinear form 
$\langle \ ,\ \rangle $ satisfies the following equation:
\begin{equation}
\label{bi-formula}
\langle X \circ Y , Z \rangle  = \langle X, Y \circ Z\rangle, 
\quad {}^{\forall} X, Y, Z \in \J. 
\end{equation}
%

For $X, Y \in \J$, the cross product $X \times Y$ is, by definition,
the element satisfying 
the following equation:
\begin{equation*}
\langle  X \times Y,  Z \rangle  = 3D(X, Y, Z),
\quad {}^{\forall} Z \in \J. 
\end{equation*}
Let $e=I_3$. Then the following equations are satisfied 
(see \cite[p.122, Lemma 5.2.1]{Springer}). 
\begin{equation}
\label{cross-formula}
{}^{\forall} X\in \J, \quad X \circ (X\times X) = \det(X)e, 
\quad e \times e = e.
\end{equation}

For any $g \in \GL(\J)$, we define $\widetilde{g} \in \GL(\J)$ by
\begin{equation*}
\langle  g(X), \widetilde{g}(Y)\rangle = \langle X, Y\rangle,
\quad {}^{\forall}  X, Y \in \J.
\end{equation*}
The following lemma is proved in 
\cite[p.180, Proposition 7.3.1] {Springer}.
\begin{lem}
\label{lem:tilde}
The map $g \mapsto \widetilde{g}$ is an automorphism of $H_1$ 
with order $2$ and for all $X,Y\in\J$, 
\begin{align*}
g(X \times Y) = \widetilde{g}(X) \times \widetilde{g}(Y), \ \,
\widetilde{g}(X \times Y)  = g(X) \times g(Y).
\end{align*}
\end{lem}

\begin{cor}
\label{cor:gtilde}
If $g\in G_1$ and $X,Y,Z,W\in\J$, then 
%
%
%
\begin{equation*}
g((X\times Y)\times (Z\times W))
= c(g)^{-1}
\left(g(X) \times g(Y)\right) \times \left(g(Z) \times g(W)\right). 
\end{equation*}
%
%
\end{cor}
\begin{proof}
There exists $t \in \overline k$ 
such that
$t^3 = c(g_1)$.
Then $t^{-1} g  \in H_{1\, \overline k}$. So 
\begin{align*}
g\left((X \times Y)\times (Z \times W)\right)  
&=
t(t^{-1}g)\left((X \times Y)\times (Z \times W)\right)\\
&=
t\left(t^{-1}g(X) \times t^{-1}g(Y)\right) 
\times \left(t^{-1}g(Z) \times t^{-1}g(W)\right)\\
&= c(g)^{-1}
\left(g(X) \times g(Y)\right) \times \left(g(Z) \times g(W)\right).
\end{align*}
\end{proof}


We now define an equivariant map 
\begin{equation*}
\text{S}:V\ni x \mapsto \mathrm{S}_x\in \text{Hom}_k(\J\otimes \J,\J)
\end{equation*}
such that $\text{S}_w(X,Y)=X\circ Y$ 
for all $X,Y\in \J$.

Let 
\begin{equation}
\label{eq:xy-defn}
X = 
\begin{pmatrix}
s_1 & x_3 & \overline {x_2} \\
\overline {x_3} & s_2 & x_1 \\
x_2 & \overline {x_1} & s_3 
\end{pmatrix},\;
Y = 
\begin{pmatrix}
t_1 & y_3 & \overline {y_2} \\
\overline {y_3} & t_2 & y_1 \\
y_2 & \overline {y_1} & t_3 
\end{pmatrix}
\end{equation}
where $s_i,t_i\in k,x_i,y_i\in \Oct$
for $i=1,2,3$. 
By computation, $X\circ Y$ is 
the following matrix:
\begin{equation}
\label{eq:XcircY-defn}
\frac 12
\begin{pmatrix}
2s_1t_1+\text{tr}(x_3\overline {y_3}+x_2\overline {y_2})
& 
\begin{matrix}
s_1y_3+t_2x_3+\overline {x_2}\, \overline {y_1} \\
+ t_1x_3+s_2y_3+\overline {y_2}\, \overline {x_1} 
\end{matrix}
& 
\begin{matrix}
s_1\overline {y_2}+x_3y_1+t_3\overline {x_2} \\ 
+t_1\overline {x_2}+y_3x_1+s_3\overline {y_2} 
\end{matrix} \\[10pt]
\begin{matrix}
s_1\overline {y_3}+t_2\overline {x_3}+ y_1x_2 \\
+ t_1\overline {x_3}+s_2\overline {y_3}+ x_1y_2 
\end{matrix}
& 2s_2t_2 + \text{tr}(x_3\overline {y_3}+x_1\overline {y_1})
& 
\begin{matrix}
\overline {x_3}\, \overline {y_2}+s_2y_1+t_3x_1 \\
+ \overline {y_3}\, \overline {x_2}+t_2x_1+s_3y_1  
\end{matrix} \\[10pt]
\begin{matrix}
s_1y_2+\overline {y_1}\, \overline {x_3}+t_3x_2 \\ 
+t_1x_2+\overline {x_1}\, \overline {y_3}+s_3y_2 
\end{matrix}
&
\begin{matrix}
y_2x_3+s_2\overline {y_1}+t_3\overline {x_1} \\
+ x_2y_3+t_2\overline {x_1}+s_3\overline {y_1}  
\end{matrix}
& 2s_3t_3+ \text{tr}(x_2\overline {y_2}+x_1\overline {y_1})
\end{pmatrix}. 
\end{equation}
Note that 
\begin{math}
\text{tr}(\overline {x_2} y_2)=\text{tr}(y_2\overline {x_2})
=\text{tr}(x_2 \overline {y_2}), \;\text{etc.} 
\end{math}
In particular, if $Y$ is diagonal, then 
\begin{equation}
\label{eq:XcircYdiag-defn}
X\circ Y = 
\frac 12 
\begin{pmatrix}
2s_1t_1
& (t_1+t_2)x_3
& (t_1+t_3)\overline {x_2} \\
(t_1+t_2) \overline {x_3} 
& 2s_2t_2 
& (t_2+t_3)x_1 \\
(t_1+t_3)x_2
& (t_2+t_3)\overline {x_1}
& 2s_3t_3
\end{pmatrix}.
\end{equation}

It is known (see \cite[p.122]{Springer}) that 
\begin{equation}
\label{eq:XtimesY-defn} 
\begin{aligned}
X\times Y 
& = X\circ Y - \frac 12 \lan X,e\ran Y - \frac 12 \lan Y,e\ran X 
- \frac 12 \lan X,Y\ran e
+ \frac 12 \lan X,e\ran \lan Y,e\ran e \\
& = X\circ Y - \frac 12 \text{Tr}(X) Y 
- \frac 12 \text{Tr}(Y) X 
- \frac 12 \text{Tr}(X\circ Y) e
+ \frac 12 \text{Tr}(X)\text{Tr}(Y) e. 
\end{aligned}
\end{equation}

For $a,b\in \J$ and $i_1\ccd i_8=1,2$, 
we define 
\begin{equation}
\label{eq:i1toi8}
[[i_1,i_2|i_3,i_4|i_5,i_6|i_7,i_8]]_{(a,b)} 
= v_1\otimes \cdots \otimes v_8
\end{equation}
where $v_j=a$ (resp. $v_j=b$) if $i_j=1$ (resp. $i_j=2$).  
For example, 
\begin{equation*}
[[1,2|2,1|2,1|1,2]]_{(a,b)}
= a\otimes b\otimes b\otimes a\otimes b\otimes a\otimes 
a\otimes b.  
\end{equation*}
We only consider 
\begin{math}
[[i_1,i_2|i_3,i_4|i_5,i_6|i_7,i_8]]_{(a,b)}  
\end{math}
such that $\{i_{2j-1},i_{2j}\}=\{1,2\}$ ($j=1\ccd 4$). 

Let 
\begin{equation*}
t_4(a,b) = \overbrace{(a\otimes b-b\otimes a) \otimes \cdots \otimes 
(a\otimes b-b\otimes a)}^4 \in \J^{8\otimes}.
\end{equation*}
Note that $t_4(a,b)$ depends only on $a\wedge b\in \wedge^2 \J$.  
By expanding all the terms,
\begin{align*}
t_4(a,b) 
& = [[1,2|1,2|1,2|1,2]]_{(a,b)} - [[1,2|1,2|1,2|2,1]_{(a,b)} \\
& \quad - [[1,2|1,2|2,1|1,2]]_{(a,b)} + [[1,2|1,2|2,1|2,1]]_{(a,b)} \\
& \quad - [[1,2|2,1|1,2|1,2]]_{(a,b)} + [[1,2|2,1|1,2|2,1]]_{(a,b)} \\
& \quad + [[1,2|2,1|2,1|1,2]]_{(a,b)} - [[1,2|2,1|2,1|2,1]]_{(a,b)} \\
& \quad - [[2,1|1,2|1,2|1,2]]_{(a,b)} + [[2,1|1,2|1,2|2,1]]_{(a,b)} \\ 
& \quad + [[2,1|1,2|2,1|1,2]]_{(a,b)} - [[2,1|1,2|2,1|2,1]]_{(a,b)} \\
& \quad + [[2,1|2,1|1,2|1,2]]_{(a,b)} - [[2,1|2,1|1,2|2,1]]_{(a,b)} \\
& \quad - [[2,1|2,1|2,1|1,2]]_{(a,b)} + [[2,1|2,1|2,1|2,1]]_{(a,b)}.  
\end{align*}

For $v_1\ccd v_8,X,Y\in \J$, we put
\begin{align*}
\Psi_1(v_1\otimes \cdots \otimes v_8)(X,Y) 
& = D(v_2,v_5,v_7)D(v_4,v_6,v_8)
(v_1\times v_3)\times  (X\times Y), \\
\Psi_2(v_1\otimes \cdots \otimes v_8)(X,Y) 
& = D(v_2,v_5,v_7)D(v_4,v_6,v_8) 
(D(v_1,v_3,X)Y+D(v_1,v_3,Y)X).
\end{align*}
Then $\Psi_i$ induces a $k$-linear map 
$\J^{8\otimes}\to \text{Hom}_k(\J\otimes \J,\J)$ 
for $i=1,2$. 
We define 
\begin{equation}
\label{eq:Phi1-defn}
\Phi_{i,(a,b)}(X,Y)
= \Psi_i(t_4(a,b))(X,Y)  
\end{equation}
for $i=1,2$. Then 
\begin{equation*}
\begin{array}{rccccc}
\Phi_i: & V & \to & (\wedge \J^2)^{4\otimes} & \to &  \text{Hom}_k(\J\otimes \J,\J)   \\
& \rotatebox{90}{$\in $} 	&& \rotatebox{90}{$\in$} 
&& \rotatebox{90}{$\in$} \\
& (a,b) & \mapsto & t_4(a,b) & \mapsto & \Phi_{i,(a,b)}  
\end{array}
\end{equation*}
is a $k$-linear map 
defined by degree $8$ polynomials of $(a,b)$.

\begin{lem}
\label{lem:equivariant-1}
$\Phi_{i,gx}(g_1X,g_1Y)=c(g_1)^3\det(g_2)^4g_1(\Phi_{i,x}(X,Y))$ 
for $i=1,2$ and all $g=(g_1,g_2)\in G$, $X,Y\in\J$. 
\end{lem}
\begin{proof}
If $x=(x_1,x_2)$, then $t_4(x_1,x_2)$ 
depends only on $x_1\wedge x_2$. If 
$g_2\in\GL(2)$ and $y=(y_1,y_2)=gx$,  
then $y_1\wedge y_2 = (\det(g_2))x_1\wedge x_2$. 
Since $x_1\otimes x_2-x_2\otimes x_1$ 
can be identified with $x_1\wedge x_2$ 
and it appears in the definition of 
$t_4(x_1,x_2)$ four times, 
$t_4(y_1,y_2)=(\det g_2)^4 t_4(x_1,x_2)$. 
So we may assume that $g_2=1$ and only consider
the action of $G_1$. 

Suppose that $v_1\otimes \cdots \otimes v_8$ 
is a term which appears in the expansion of 
$t_4(x_1,x_2)$. 
If $g_1\in G_1$ and $x$ is replaced by $g_1x=(g_1x_1,g_1x_2)$, 
then terms which appear in the expansion of 
$t_4(g_1x_1,g_1x_2)$ are in the form  
$g_1v_1\otimes \cdots \otimes g_1v_8$. 
By Corollary \ref{cor:gtilde}, 
\begin{align*}
& D(g_1 v_2,g_1 v_5,g_1v_7)D(g_1v_4,g_1v_6,g_1v_8)
(g_1v_1\times g_1v_3)\times  (g_1X\times g_1Y), \\ 
& \quad = c(g_1)^3 D(v_2,v_5,v_7)D(v_4,v_6,v_8)
g_1((v_1\times v_3)\times  (X\times Y)). 
\end{align*}
Also 
\begin{align*}
& D(g_1v_2,g_1v_5,g_1v_7)D(g_1v_4,g_1v_6,g_1v_8)
(D(g_1v_1,g_1v_3,g_1X)g_1Y+D(g_1v_1,g_1v_3,g_1Y)g_1X) \\
& \quad = c(g_1)^3D(v_2,v_5,v_7)D(v_4,v_6,v_8) 
g_1(D(v_1,v_3,X)Y+D(v_1,v_3,Y)X).
\end{align*}
Therefore, 
$\Phi_{i,g_1x}(g_1X,g_1Y)=c(g_1)^3g_1(\Phi_{i,x}(X,Y))$
for $i=1,2$. 
\end{proof}

We first evaluate $\Phi_{1,w}$ (see (\ref{eq:w-defn})). 
\begin{prop}
\label{prop:Phi1evaluated} 
For all $X,Y\in\J$, 
\upshape
\begin{math}
-18\Phi_{1,w}(X,Y) 
= X\circ Y - \frac 12 \text{Tr}(Y)X - \frac 12 \text{Tr}(X)Y. 
\end{math}
\end{prop}
\begin{proof}
Suppose that $(a,b)=w=(w_1,w_2)$ in the following. 
If we expand $t_4(a,b)$, then the coefficient of 
\begin{math}
[[i_1,i_2|i_3,i_4|i_5,i_6|i_7,i_8]]_{(a,b)}  
\end{math}
is $1$ (resp. $-1$) if the number of $j=1\ccd 4$ 
such that $(i_{2j-1},i_{2j})=(2,1)$ is even (resp. odd). 
We list 
the sign in $t_4(a,b)$, $D(v_2,v_5,v_7)D(v_4,v_6,v_8)$
and $v_1\times v_3$ 
in the following table assuming that
\begin{math}
[[i_1,i_2|i_3,i_4|i_5,i_6|i_7,i_8]]_{(a,b)}  
\end{math}
is in the form (\ref{eq:i1toi8}).


\begin{center}
\begin{tabular}{|c|c|c|}
\hline
(1)
& $D(v_2,v_5,v_7)D(v_4,v_6,v_8)$
& $v_1\times v_3$ \\
\hline
$+ [[1,2|1,2|1,2|1,2]]$  
& $D(b,a,a) D(b,b,b)$  
& $a\times a$ \\
\hline
$- [[1,2|1,2|1,2|2,1]]$ 
& $D(b,a,b) D(b,b,a)$  
& $a\times a$ \\
\hline
$- [[1,2|1,2|2,1|1,2]]$ 
& $D(b,b,a) D(b,a,b)$  
& $a\times a$ \\
\hline
$+ [[1,2|1,2|2,1|2,1]]$  
& $D(b,b,b) D(b,a,a)$  
& $a\times a$ \\
\hline
$- [[1,2|2,1|1,2|1,2]]$  
& $D(b,a,a) D(a,b,b)$  
& $a\times b$ \\
\hline
$+ [[1,2|2,1|1,2|2,1]]$ 
& $D(b,a,b) D(a,b,a)$  
& $a\times b$ \\
\hline
$+ [[1,2|2,1|2,1|1,2]]$ 
& $D(b,b,a) D(a,a,b)$  
& $a\times b$ \\
\hline
$- [[1,2|2,1|2,1|2,1]]$ 
& $D(b,b,b) D(a,a,a)$  
& $a\times b$ \\
\hline
$- [[2,1|1,2|1,2|1,2]]$ 
& $D(a,a,a) D(b,b,b)$  
& $b\times a$ \\
\hline
$+ [[2,1|1,2|1,2|2,1]]$ 
& $D(a,a,b) D(b,b,a)$  
& $b\times a$ \\
\hline
$+ [[2,1|1,2|2,1|1,2]]$ 
& $D(a,b,a) D(b,a,b)$  
& $b\times a$ \\
\hline
$- [[2,1|1,2|2,1|2,1]]$ 
& $D(a,b,b) D(b,a,a)$  
& $b\times a$ \\
\hline
$+ [[2,1|2,1|1,2|1,2]]$ 
& $D(a,a,a) D(a,b,b)$  
& $b\times b$ \\
\hline
$- [[2,1|2,1|1,2|2,1]]$ 
& $D(a,a,b) D(a,b,a)$  
& $b\times b$ \\
\hline
$- [[2,1|2,1|2,1|1,2]]$ 
& $D(a,b,a) D(a,a,b)$  
& $b\times b$ \\
\hline
$+ [[2,1|2,1|2,1|2,1]]$ 
& $D(a,b,b) D(a,a,a)$  
& $b\times b$ \\
\hline
\end{tabular}
\end{center}

\noindent
((1) is $[[i_1,i_2|i_3,i_4|i_5,i_6|i_7,i_8]]_{(a,b)}$ with its sign in 
$t_4(a,b)$.)

Note that if $(a,b)=w$, then 
\begin{equation*}
D(a,a,a) = D(b,b,b)=0,\;
D(a,a,b) = \frac 13,\;
D(a,b,b) = -\frac 13.  
\end{equation*}
So we can ignore terms with the second column including 
either $D(a,a,a)$ or $D(b,b,b)$.  
Removing these terms, we obtain the following table. 

\begin{center}
\begin{tabular}{|c|c|c|}
\hline
(1)
& $9D(v_2,v_5,v_7)D(v_4,v_6,v_8)$
& $v_1\times  v_3$ \\
\hline
$- [[1,2|1,2|1,2|2,1]]$ 
& $1$  
& $a\times a$ \\
\hline
$- [[1,2|1,2|2,1|1,2]]$ 
& $1$  
& $a\times a$ \\
\hline
$- [[1,2|2,1|1,2|1,2]]$  
& $-1$  
& $a\times b$ \\
\hline
$+ [[1,2|2,1|1,2|2,1]]$ 
& $-1$  
& $a\times b$ \\
\hline
$+ [[1,2|2,1|2,1|1,2]]$ 
& $-1$  
& $a\times b$ \\
\hline
$+ [[2,1|1,2|1,2|2,1]]$ 
& $-1$  
& $b\times a$ \\
\hline
$+ [[2,1|1,2|2,1|1,2]]$ 
& $-1$  
& $b\times a$ \\
\hline
$- [[2,1|1,2|2,1|2,1]]$ 
& $-1$  
& $b\times a$ \\
\hline
$- [[2,1|2,1|1,2|2,1]]$ 
& $1$  
& $b\times b$ \\
\hline
$- [[2,1|2,1|2,1|1,2]]$ 
& $1$  
& $b\times b$ \\
\hline
\end{tabular}
\end{center}

By the above table, 
\begin{align*}
9\Phi_{1,w}(X,Y)
& =  -(2a\times a+ a\times b+a\times b+2b\times b)\times (X\times Y).
\end{align*}
If we put 
\begin{math}
c = a+b = 
\text{diag}(1,0,-1),
\end{math}
then 
\begin{align*}
9\Phi_{1,w}(X,Y)
& =  -(a\times a+b\times b+c\times c)\times (X\times Y). 
\end{align*}

By (\ref{eq:XcircYdiag-defn}) and (\ref{eq:XtimesY-defn}) 
\begin{align*}
a\times a & = \text{diag}(0,0,-1), \;
b\times b  = \text{diag}(-1,0,0), \;
c\times c  = \text{diag}(0,-1,0).
\end{align*}
Therefore, 
\begin{math}
9\Phi_{1,w}(X,Y)
= e\times (X\times Y).  
\end{math}
Since $e$ is the unit element of $\J$, 
replacing $X,Y$ in (\ref{eq:XtimesY-defn})
by $e,X\times Y$, we obtain 
\begin{align*}
9\Phi_{1,w}(X,Y) 
& = X\times Y - \frac 12 \text{Tr}(e) X\times Y 
- \frac 12 \text{Tr}(X\times Y) e 
- \frac 12 \text{Tr}(X\times Y) e 
+ \frac 12 \text{Tr}(e)\text{Tr}(X\times Y) e. \\
& = -\frac 12  X\times Y 
+ \frac 12 \text{Tr}(X\times Y) e. 
\end{align*}
By (\ref{eq:XtimesY-defn}), 
\begin{align*}
\text{Tr}(X\times Y)
& = \text{Tr}(X\circ Y)-\text{Tr}(X) \text{Tr}(Y)
-\frac 32\text{Tr}(X\circ Y) + \frac 32\text{Tr}(X) \text{Tr}(Y) \\
& = -\frac 12\text{Tr}(X\circ Y)
+ \frac 12\text{Tr}(X) \text{Tr}(Y). 
\end{align*}
Therefore, again by (\ref{eq:XtimesY-defn}), 
\begin{equation}
\label{eq:Phi1evaluated} 
9\Phi_{1,w}(X,Y) 
= -\frac 12  X\circ Y + \frac 14 \text{Tr}(Y)X + \frac 14 \text{Tr}(X)Y. 
\end{equation}
Multiplying $-2$, we obtain the proposition. 
\end{proof}

\begin{prop}
\label{prop:Phi2-evaluate}
For all $X,Y\in \J$, 
\upshape
\begin{math}
3\Phi_{2,w}(X,Y) 
= \text{Tr}(Y)X+\text{Tr}(X)Y. 
\end{math} 
\end{prop}
\begin{proof}
By very similar comutations as in the case of  
$\Phi_{1,w}(X,Y)$, 
we obtain 
\begin{equation}
\label{eq:Phi2evaluated} 
\begin{aligned}
3\Phi_{2,w}(X,Y)
& = -3(D(a,a,X) + D(b,b,X) + D(c,c,X))Y \\ 
& \quad -3(D(a,a,Y) + D(b,b,Y) + D(c,c,Y))X \\
& = (s_3+s_1+s_2)Y + (t_3+t_1+t_2)X \\
& = \text{Tr}(Y)X+\text{Tr}(X)Y.
\end{aligned}
\end{equation}
\end{proof}

For $(a,b)\in V_k$, we define 
\begin{equation}
\label{eq:C1-defn}
{\mathrm S}_{(a,b)}(X,Y) = -18\Phi_{1,(a,b)}(X,Y) + \frac 32 \Phi_{2,(a,b)}(X,Y). 
\end{equation}

The following proposition follows from 
Lemma \ref{lem:equivariant-1} and 
(\ref{eq:Phi1evaluated}), (\ref{eq:Phi2evaluated}). 

\begin{prop}
\label{prop:equivariantS}
\begin{itemize}
\item[(1)] 
If $x\in V_k$, $g=(g_1,g_2)\in G_k$ and 
$X,Y\in \J$, then 
\begin{equation*}
{\mathrm S}_{gx}(g_1X,g_1Y)=c(g_1)^3\det(g_2)^4g_1({\mathrm S}_x(X,Y)). 
\end{equation*}
\item[(2)] 
${\mathrm S}_w(X,Y) = X\circ Y$. 
\end{itemize} 
\end{prop}

Let $\Delta(x)$ be the relaaive invariant polynomial of
degree $12$ and $\Delta(w)=1$ which we defined in Introduction.
For $x\in V^{\mathrm{ss}}_k$, we define a $k$-algebra 
structure $\circ_x$ on $\J$ by 
\begin{equation}
\label{eq:x-algebra}
X\circ_x Y \stackrel {\rm{def}} {=} \Delta(x)^{-1}{\mathrm S}_x(X,Y).
\end{equation}
We denote this $k$-algebra on the underlying vector space $\J$ 
by $\Mj_x$. Proposition \ref{prop:equivariantS} (2) implies that 
$\Mj_w$ is isomorphic to $\J$ (the original Jordan algebra structure).

For $g = (g_1, g_2)\in G$, 
we define an element of $G_1$ by
\begin{equation}
\label{eq:mu-g-defn}
\mu_g = c(g_1) \det(g_2)^2 g_1 \in G_1.
\end{equation}
Note that the map $g \to \mu_g$ is a homomorphism.

\begin{thm}
 \label{thm:main-th1}
If $x,y\in V^{\mathrm{ss}}_k$, $g\in G_{\ks}$ and
$y=gx$, then $\mu_g:\Mj_{x\,\ks}\to \Mj_{y\,\ks}$
is an isomorphism of $k$-algebras
\end{thm}
\begin{proof}
Let $X,Y\in \J_{\ks}=\Mj_{x\,\ks}$. Then 
\begin{align*}
\mu_g(X) \circ_y \mu_g(Y) 
& = \Delta(gx)^{-1} {\mathrm S}_{gx}
(c(g_1)\det(g_2)^2g_1(X),c(g_1)\det(g_2)^2g_1(Y)) \\
& = \Delta(gx)^{-1} c(g_1)^2\det(g_2)^4{\mathrm S}_{gx}
(g_1(X),g_1(Y)) \\
& = \Delta(gx)^{-1} c(g_1)^2\det(g_2)^4
c(g_1)^3 \det(g_2)^4 g_1({\mathrm S}_x(X,Y)) 
\;\text{by Proposition \ref{prop:equivariantS}} \\
& = \Delta(gx)^{-1}\Delta(x) c(g_1)^4\det(g_2)^6
c(g_1)\det(g_2)^2 g_1(\Delta(x)^{-1}{\mathrm S}_x(X,Y)) \\
& = c(g_1)\det(g_2)^2 g_1(X\circ_x Y)=\mu_g(X\circ_x Y). 
\end{align*}
Therefore, $\mu_g$ is a homomorphism. Since $\mu_g$
is obviously bijecive, it is an isomorphism. 
\end{proof}

The above theorem implies that if $g_x=(g_1,g_2)\in G_{\ks}$
and $x=g_xw\in V^{\mathrm{ss}}_k$, then 
$\mu_g^{-1}(\Mj_x)\subset \J_{\ks}=\Mj_{w\,\ks}$ is a $k$-form of 
$\J$. If $(\Mj,\an)\in{\mathrm{JIC}}(k)$ is the pair 
corresponding to $x$ (see \cite[Section 4,5]{kato-yukie-jordan}), 
$\Mj$ was characterized in the same manner
(see \cite[(4.9)]{kato-yukie-jordan}).
Therefore, Theorem \ref{thm:main-theorem-equiv} follows.


\section{Equivariant map II}
\label{sec:equivariant-mapII}

In this section, we prove Theorem \ref{thm:main-theorem-equiv2}. 

Let 
\begin{equation}
\mathrm{T}_a(X,Y,Z)
= 27D(a,a,X)D(a,a,Y)D(a,a,Z) 
- 24D(a,a,a)D(a\times X,a\times Y,a\times Z)    
\end{equation}
for $a,X,Y,Z\in \J$. 
Then the map 
\begin{equation*}
\mathrm{T}: \J \ni a \mapsto 
\mathrm{T}_a   \in 
{\mathrm{Hom}}_k(\J\otimes \J\otimes \J,k) 
\end{equation*}
is $k$-linear.

\begin{lem}
If $a\in \J$, $g\in G_1$ and 
$X,Y,Z\in \J$, then 
\begin{equation*}
{\mathrm T}_{ga}(gX,gY,gZ)=c(g)^3g({\mathrm T}_a(X,Y,Z)). 
\end{equation*}
\end{lem}
\begin{proof}
There exists $t\in \overline k$ such that 
$t^3=c(g)$. We put $g_1 = t^{-1}g$. 
Then $g=tg_1$ and $g_1\in H_{1\,\overline k}$. 
So, 
\begin{align*}
D(ga\times gX,ga\times gY,ga\times gZ)    
& = t^6D(g_1a\times g_1X,g_1a\times g_1Y,g_1a\times g_1Z) \\   
& = t^6D(\widetilde g_1(a\times X),
\widetilde g_1(a\times Y),
\widetilde g_1(a\times Z)) \\
& = t^6D(a\times X,a\times Y,a\times Z) \\
& = c(g)^2D(a\times X,a\times Y,a\times Z).
\end{align*}
Since 
\begin{math}
D(ga,ga,gX)= c(g)D(a,a,X),  
\end{math}
etc., the lemma follows. 
\end{proof}

The following proposition plays a crucial role 
in proving Theorem \ref{thm:main-theorem-equiv2}. 

\begin{prop}
\label{prop:equiv-fromJ}
$\mathrm{T}_e(X,Y,Z)=\mathrm{Tr}((X\circ Y)\circ Z)$ 
($= \mathrm{Tr}(X\circ (Y\circ Z))$) 
for all $X,Y,Z\in \J$.  
\end{prop}
\begin{proof}
Since $D(e,e,e)=1$, 
\begin{equation*}
\mathrm{T}_e(X,Y,Z)
=  27D(e,e,X)D(e,e,Y)D(e,e,Z)
-24D(e\times X,e\times Y,e\times Z). 
\end{equation*}
Let 
\begin{equation*}
X = 
\begin{pmatrix}
s_1 & x_3 & \overline {x_2} \\
\overline {x_3} & s_2 & x_1 \\
x_2 & \overline {x_1} & s_3 
\end{pmatrix},\;
Y = 
\begin{pmatrix}
t_1 & y_3 & \overline {y_2} \\
\overline {y_3} & t_2 & y_1 \\
y_2 & \overline {y_1} & t_3 
\end{pmatrix},\;
Z = 
\begin{pmatrix}
u_1 & z_3 & \overline {z_2} \\
\overline {z_3} & u_2 & z_1 \\
z_2 & \overline {z_1} & u_3 
\end{pmatrix}  \in\J.
\end{equation*}

Note that 
\begin{align*}
6D(X,Y,Z)
& = \sum_{\{i,j,k\}=\{1,2,3\}} s_i t_j u_k 
+ \sum_{\{i,j,k\}=\{1,2,3\}} \text{tr} (x_i y_j z_k) \\
& \quad - \sum_i s_i\text{tr}(y_i\overline {z_i})
- \sum_i t_i\text{tr}(x_i\overline {z_i})
- \sum_i u_i\text{tr}(x_i\overline {y_i}). 
\end{align*}
By (\ref{eq:XcircYdiag-defn}) and 
(\ref{eq:XtimesY-defn}), 
\begin{align*}
e\times X 
& = -\frac 12 
\begin{pmatrix}
-(s_2+s_3) & x_3 & \overline {x_2} \\
\overline {x_3} & -(s_1+s_3) & x_1 \\
x_2 & \overline {x_1} & -(s_1+s_2)  
\end{pmatrix}
\end{align*}

Therefore, 
\begin{align*}
48D(e\times X,e\times Y,e\times Z) 
& = [(s_2+s_3)(t_1+t_3)(u_1+u_2)
+ (s_2+s_3)(t_1+t_2)(u_1+u_3) \\
& \quad + (s_1+s_3)(t_2+t_3)(u_1+u_2)
+ (s_1+s_3)(t_1+t_2)(u_2+u_3) \\
& \quad + (s_1+s_2)(t_2+t_3)(u_1+u_3)
+ (s_1+s_2)(t_1+t_3)(u_2+u_3)] \\
& \quad - \sum_{\{i,j,k\}=\{1,2,3\}}\text{tr}(x_iy_jz_k) \\
 & \quad - (s_2+s_3)\text{tr}(y_1\overline {z_1})
- (s_1+s_3)\text{tr}(y_2\overline {z_2})
- (s_1+s_2)\text{tr}(y_3\overline {z_3}) \\
& \quad - (t_2+t_3)\text{tr}(x_1\overline {z_1})
+ (t_1+t_3)\text{tr}(x_2\overline {z_2})
+ (t_1+t_2)\text{tr}(x_3\overline {z_3}) \\
& \quad - (u_2+u_3)\text{tr}(x_1\overline {y_1})
+ (u_1+u_3)\text{tr}(x_2\overline {y_2})
+ (u_1+u_2)\text{tr}(x_3\overline {y_3}) \\
& = 2\sum_{\{i,j,k\}=\{1,2,3\}} s_it_ju_k
+ 2(s_1t_1u_2+\cdots) 
- \sum_{\{i,j,k\}=\{1,2,3\}}\text{tr}(x_iy_jz_k) \\
& \quad - \sum_{i\not=j} s_i\text{tr}(y_j\overline {z_j})
- \sum_{i\not=j} t_i\text{tr}(x_j\overline {z_j})
- \sum_{i\not=j} u_i\text{tr}(x_j\overline {y_j}). 
\end{align*}

Since 
\begin{math}
6D(e,e,X) = 2(s_1+s_2+s_3), 
\end{math}
we have 
\begin{align*}
27D(e,e,X)D(e,e,Y)D(e,e,Z) 
& = (s_1+s_2+s_3)
(t_1+t_2+t_3) (u_1+u_2+u_3) \\
& = \sum_i s_it_iu_i  
+ (s_1t_1u_2+\cdots)
+ \sum_{\{i,j,k\}=\{1,2,3\}} s_it_ju_k.
\end{align*}
Therefore, 
\begin{equation}
\label{eq:Te=Tri}
\begin{aligned}
\mathrm{T}_e(X,Y,Z)
& = \sum_i s_it_iu_i 
+ \frac 12 \sum_{\{i,j,k\}=\{1,2,3\}}\text{tr}(x_iy_jz_k) \\
& \quad + \frac 12 \sum_{i\not=j} s_i\text{tr}(y_j\overline {z_j})
+ \frac 12 \sum_{i\not=j} t_i\text{tr}(x_j\overline {z_j})
+ \frac 12 \sum_{i\not=j} u_i\text{tr}(x_j\overline {y_j}).
\end{aligned}
\end{equation}
%

%
%
%
%
%
%

Replacing $X,Y$ in (\ref{eq:XcircY-defn}) by $X\circ Y,Z$
respetively and taking the sum of diagonal entries, 
we can express 
$4\mathrm{Tr}((X\circ Y)\circ Z)$ in the following manner:
\begin{align*}
& 2(2s_1t_1+\text{tr}(x_3\overline y_3+\overline x_2y_2))u_1 \\
& + \text{tr}(
(s_1y_3+x_3t_2+\overline {x_2}\, \overline {y_1} 
+ t_1x_3+y_3s_2+\overline {y_2}\, \overline {x_1})\overline {z_3}
+(s_1\overline {y_2}+x_3y_1+\overline {x_2}t_3  
+t_1\overline {x_2}+y_3x_1+\overline {y_2}s_3)z_2) \\
& + 2(2s_2t_2+\text{tr}(x_1\overline {y_1}+\overline {x_3}y_3))u_2 \\
& + \text{tr}(
(\overline {x_3}\, \overline {y_2}+s_2y_1+x_1t_3
+\overline {y_3}\, \overline {x_2}+t_2x_1+y_1s_3)\overline {z_1}
+(\overline {x_3}t_1+s_2\overline {y_3}+x_1y_2 
+ \overline {y_3}s_1+t_2\overline {x_3}+y_1x_2)z_3) \\
& + 2(2s_3t_3+\text{tr}(x_2\overline {y_2}+\overline {x_1}y_1))u_3 \\ 
& + \text{tr}(
(x_2t_1+\overline {x_1}\, \overline {y_3}+s_3y_2 
+ y_2s_1+\overline {y_1}\, \overline {x_3}+t_3x_2)\overline {z_2}
+(x_2y_3+\overline {x_1}t_2+s_3\overline {y_1} 
+ y_2x_3+\overline {y_1}t_2+s_3\overline {x_1})z_1). 
\end{align*}
This coincides with $4$ times (\ref{eq:Te=Tri}). 
\end{proof}

The construction of the isotope 
defined for $a\in\J$ ($\det(a)\not=0$) 
is given in \cite[p.155]{Springer}. 
Let 
\begin{equation}
\label{eq:bilinear-product}
\begin{aligned}
Q_a(X,Y)
& =  -6\det(a)D(X, Y, a)  
+ 9D(X,a,a)D(Y,a,a), \\
\Phi_a(X,Y)
& =  4\det(a)^3 (X \times a) \times (Y \times a) 
+ \frac 12
(\det(a)^2Q_a(X,Y)
-Q_a(X,a)Q_a(Y,a))a. 
\end{aligned}
\end{equation}
Then the product structure and the associated  
bilinear form of the isotope corresponding to $a\in\J$ is 
given by 
\begin{align*}
& \J^2 \ni (X,Y)\mapsto X\circ_a Y
\stackrel {\rm{def}} {=} \det(a)^{-4}\Phi_a(X,Y)\in\J, \\
& \J^2 \ni (X,Y)\mapsto \langle X,Y\rangle_a
\stackrel {\rm{def}} {=} \det(a)^{-2}Q_a(X,Y)\in k.
\end{align*}
(see \cite[p.147, Proposition 5.6.2]{Springer}, 
\cite[p.153, Proposition 5.8.2]{Springer}, 
\cite[p.155, Proposition 5.9.2]{Springer}).

Here we provide an alternative way to define the product
structure on $\J_a$.  

\begin{mydef}
\label{defn:Ja-structure-defn}
Suppose that $a\in \J,\det(a)\not=0$. 
For $X,Y\in \J$, we define $X\circ_a Y\in\J$ to be the element
such that $Q_a(X\circ_a Y,Z)=\det(a)^{-1}{\mathrm T}_a(X,Y,Z)$
for all $Z\in\J$.  
\end{mydef}

Since $Q_a$ is a non-degenerate bilinear form, the definition of 
$X\circ_a Y$ is well-defined. 

\begin{thm}
\label{thm:main-th2}
If $a\in\J$, $\det(a)\not=0$, 
$g_a\in GE_{6\,\ks}$, $a=g_ae\in \J$, 
then $\J\otimes \ks \ni X\mapsto gX\in \J_a\otimes \ks$ 
induces an isomorphism of $k$-algebras 
$\J\otimes \ks\cong \J_a\otimes \ks$.  
\end{thm}
\begin{proof}
Suppose that $X,Y,Z\in \J\otimes \ks$. 
Since $\det(a)=c(g)\det(e)=c(g)$, 
\begin{align*}
Q_a(g(X\circ Y),gZ) 
& = c(g)^2 Q_e(X\circ Y,Z) 
= c(g)^2 {\mathrm T}_e(X,Y,Z)
= c(g)^{-1}{\mathrm T}_a(gX,gY,gZ) \\
& = \det(a)^{-1}{\mathrm T}_a(gX,gY,gZ) = Q_a(gX\circ_a gY,gZ). 
\end{align*}
Since this holds for all $Z\in\J\otimes\ks$, 
$g(X\circ Y)=gX\circ_a gY$. Therefore, 
$g$ induces an isomorphism of $k$-algebras
$\J\otimes \ks\cong \J_a\otimes \ks$. 
\end{proof}

In the situation of Theorem \ref{thm:main-th2}, 
$\J_a$ can be identified with the $k$-form
$g^{-1}(\J_a)\subset \J_{\ks}$.
Therefore, this $\J_a$ concides
with the isotope $\J_a$ constructed in \cite{Springer}.

Let $m:V\to\J$ be the equivariant map in \cite[Section 5]{kato-yukie-jordan}.
Since $m$ is defined by homogeneous polynomials of degree $4$, 
$Q_{m(x)},{\mathrm T}_{m(x)}$ are defined by homogeneous polynomials of 
degrees $16,24$ respectively. Can we construct lower degree 
equivariant maps from $V$ to the space of bilinear forms and 
trilinear forms on $\J$ to give an alternative way 
to define the product sutrucure on $\Mj_x$ 
(see (\ref{eq:x-algebra}))? It does not seem so. 
For example, if there were such a quadratic form $Q_x$
depending on $x$, the degree must be $16-12=4$. 
However, there seems to be only one equivariant map of degree $4$
from $\wedge^2 \J$ to the space of bilinear forms on $\J$
by the calculation of the software ``LiE'' (\cite{liesoft})
as follows.

\begin{verbatim}
 > alt_tensor(2,[1,0,0,0,0,0])
      1X[0,0,1,0,0,0]
 > sym_tensor(2,[0,0,1,0,0,0])
      1X[0,0,0,0,0,2] +1X[0,0,2,0,0,0] +1X[0,1,0,0,1,0] +1X[1,0,0,0,0,0] +
      1X[1,1,0,0,0,0] +1X[2,0,0,0,0,1]
\end{verbatim}

We can consutruct an equivariant map from 
$\wedge^2 \J$ to the space of bilinear forms on $\J$
such that $w_1\wedge w_2$ corresponds to 
the bilinear form 
\begin{equation*}
\J^2 \ni (X,X) \mapsto 
\sum_{i=1}^3s_i^2 -\sum_{i=1}^3 \|x_i\|\in k
\end{equation*}
where as 
\begin{equation*}
{\mathrm{Tr}}(X\circ X)
= \sum_{i=1}^3s_i^2 + \sum_{i=1}^3 \|x_i\|.  
\end{equation*}
Therefore, we cannot obtain $\mathrm{Tr}(X\circ X)$.

\bibliographystyle{plain}
\bibliography{ref-rep7.bib}

\end{document}